\documentclass[a4paper]{article}
\usepackage{amsfonts}
\usepackage{amsmath}
\usepackage{amssymb}
\begin{document}
\newtheorem{theorem}{Theorem}[section]
\newtheorem{lemma}[theorem]{Lemma}
\newtheorem{corollary}[theorem]{Corollary}
\newtheorem{conjecture}[theorem]{Conjecture}
\newtheorem{remark}[theorem]{Remark}
\newtheorem{definition}[theorem]{Definition}
\newtheorem{problem}[theorem]{Problem}
\newtheorem{example}[theorem]{Example}
\newtheorem{proposition}[theorem]{Proposition}
\title{{\bf Canonical singular hermitian metrics \\ on relative canonical bundles }}
\date{November 5, 2007}
\author{Hajime TSUJI}
\maketitle
\begin{abstract}
\noindent We introduce a new class of canonical  AZD's (called the supercanonical AZD's)  on the canonical bundles of smooth projective varieties with pseudoeffective canonical classes.  
We study the variation of the supercanonical AZD $\hat{h}_{can}$ under 
projective deformations and give a new proof of the invariance of plurigenera. 
This paper is a continuation of \cite{tu9}. 

\end{abstract}
\tableofcontents
\section{Introduction}
Let $X$ be a smooth projective variety and let $K_{X}$ be the canonical bundle 
of $X$.  
In algebraic geometry, the canonical ring $R(X,K_{X}) : = \oplus_{m=0}^{\infty}\Gamma (X,{\cal O}_{X}(mK_{X}))$ is one of the main object to study.   

Let $X$ be a smooth projective variety such that $K_{X}$ is pseudoeffective.
The purposes of this article  are twofold.  
The first purpose   is to  construct a singular hermitian metric $\hat{h}_{can}$ on $K_{X}$ such that 
\begin{enumerate}
\item $\hat{h}_{can}$ is uniquely determined by $X$.   
\item $\Theta_{\hat{h}_{can}}$ is semipositive in the sense of current.
\item $H^{0}(X,{\cal O}_{X}(mK_{X})\otimes {\cal I}(\hat{h}_{can}^{m}))
\simeq H^{0}(X,{\cal O}_{X}(mK_{X}))$ holds for every $m\geqq 0$,  
\end{enumerate}
where 
${\cal I}(\hat{h}_{can}^{m})$ denotes the multiplier ideal sheaf of 
$\hat{h}_{can}^{m}$ as is defined in \cite{n}. 
And the second purpose is to study the behavior of $\hat{h}_{can}$ on projective families.
We may summerize the 2nd and the 3rd conditions by introducing the following 
notion. 

\begin{definition}{\em {\bf (AZD)}}\label{azd}(\cite{tu,tu2})
Let $M$ be a compact complex manifold and let $L$ be a holomorphic line bundle
on $M$.  A singular hermitian metric $h$ on $L$ is said to be 
an analytic Zariski decomposition (AZD in short), if the followings hold.
\begin{enumerate}
\item $\Theta_{h}$ is a closed positive current. 
\item For every $m\geq 0$, the natural inclusion
\[
H^{0}(M,{\cal O}_{M}(mL)\otimes{\cal I}(h^{m}))\rightarrow
H^{0}(M,{\cal O}_{M}(mL))
\]
is an isomorphim.
\end{enumerate} 
$\square$
\begin{remark}\label{rm}
A line bundle $L$ on a projective manifold  $X$ admits an AZD, if and only if 
$L$ is pseudoeffective (\cite[Theorem 1.5]{d-p-s}). $\square$
\end{remark}
\end{definition} 

\noindent In this sense, the first purpose of this article is to construct an AZD
on $K_{X}$  depending only on $X$, 
when $K_{X}$ is pseudoeffective (by Remark \ref{rm} this is the minimum requirement for the existence of an AZD).

The main motivation to construct such a singular hermitian metric is to study the  canonical ring in terms of it.  
This is indeed possible.   For example, we obtain the invariance of plurigenera  under smooth projective deformations (cf. Corollary \ref{pg}). 
  
In fact the hermitian metric constructed here  is useful in many other contexts. Other applications and a generalization to subKLT pairs will be treated in the forthcoming papers (\cite{log}).     

I would like to express thanks to Professor Bo Berndtsson who pointed out 
an error in the previous version. 
\subsection{Canonical AZD $h_{can}$}

If we assume the stronger assumption that $X$ has nonnegative Kodaira dimension, we have already konwn how to construct a canonical AZD for $K_{X}$.
Let us review the construction in \cite{tu9}.  
\begin{theorem} (\cite{tu9})\label{canazd}
Let $X$ be a smooth projective variety with nonnegative Kodaira dimension. 
We set for every point $x\in X$
\[
K_{m}(x) :=  \sup \{\mid \sigma\mid^{\frac{2}{m}}(x) ;
\sigma \in \Gamma (X,{\cal O}_{X}(mK_{X})), \mid\!\!\int_{X}(\sigma\wedge\bar{\sigma})^{\frac{1}{m}}\!\mid = 1\}
\]
and 
\[
K_{\infty}(x) := \limsup_{m\rightarrow\infty} K_{m}(x).
\] 
Then 
\[
h_{can}:= \mbox{\em the lower envelope of}\,\,K_{\infty}^{-1}
\]
is an AZD on $K_{X}$. $\square$
\end{theorem}
\begin{remark} By the ring structure of $R(X,K_{X})$, we see that 
\[
 \limsup_{m\rightarrow\infty} K_{m}(x) = \sup_{m\geqq 1}K_{m}(x)
\]
holds. $\square$
\end{remark}
\begin{remark}
Since $h_{\infty}$ depends only on $X$,   
the volume
\[
\int_{X}h_{can}^{-1}
\]
is an invariant of $X$.  $\square$
\end{remark}
Apparently this construction is very canonical, i.e., $h_{can}$ depends only on 
the complex structure of $X$. 
We call $h_{can}$ the {\bf canonical AZD} of $K_{X}$.
But this construction works only if we know that the Kodaira dimension of $X$
is nonnegative apriori.  This is the main defect of $h_{can}$. 
For example, $h_{can}$ is useless to solve the abundance conjecture.

\subsection{Supercanonical AZD $\hat{h}_{can}$}\label{construction}

To avoid the defect of $h_{can}$ we introduce the new AZD $\hat{h}_{can}$.
  Let us use the following terminology. 

\begin{definition}\label{pe}
Let $(L,h_{L})$ be a singular hermitian line bundle on a complex manifold $X$.
$(L,h_{L})$ is said to be pseudoeffective, if the curvature current of $h_{L}$
is semipositive. $\square$
\end{definition}

Let $X$ be a smooth projective $n$-fold such that the canonical bundle $K_{X}$ 
is pseudoeffective.  Let $A$ be a sufficiently  ample line bundle 
such that for every pseudoeffective singular hermitian line bundle $(L,h_{L})$
on $X$,
${\cal O}_{X}(A + L)\otimes {\cal I}(h_{L})$ and 
 ${\cal O}_{X}(K_{X}+ A + L)\otimes {\cal I}(h_{L})$ are globally generated. 
Such an ample line bundle $A$ extists by $L^{2}$-estimates.  
Let $h_{A}$ be a a $C^{\infty}$ hermitian metric on $A$ 
with strictly positive curvature \footnote{Later we shall
also consider the case that $h_{A}$ is any $C^{\infty}$ hermitian metric 
(without positivity of curvature) or  a singular hermitian metric on $A$.}.
Let us fix a $C^{\infty}$ volume form $dV$ on $X$. 
By the $L^{2}$-extension theorem (\cite{o}) we may and do assume that 
$A$ is sufficiently ample so that for every $x\in X$ and 
for every pseudoeffective singular hermitian line bundle $(L,h_{L})$, there exists a bounded interpolation operator
\[
I_{x} : A^{2}(x,(A + L)_{x},h_{A}\cdot h_{L},\delta_{x})
\rightarrow A^{2}(X,A + L,h_{A}\cdot h_{L},dV)
\]
such that the operator norm of $I_{x}$ is bounded by 
a positive constant independent of $x$ and $(L,h_{L})$, 
where $A^{2}(X,A + L,h_{A}\cdot h_{L},dV)$ denotes the Hilbert space defined by
\[
 A^{2}(X,A + L,h_{A}\cdot h_{L},dV) := \{ \sigma\in \Gamma (X,{\cal O}_{X}(A + L)\otimes {\cal I}(h_{L}))\mid \int_{X}\mid\sigma\mid^{2}\cdot h_{A}\cdot h_{L}\cdot dV < + \infty\} 
\]
with the $L^{2}$ inner product 
\[
(\sigma ,\sigma^{\prime}):= \int_{X}\sigma\cdot\bar{\sigma}^{\prime}\cdot h_{A}\cdot h_{L}\cdot dV 
\]
and  $A^{2}(x,(A + L)_{x},h_{A}\cdot h_{L},\delta_{x})$ is defined 
similarly, where $\delta_{x}$ is the Dirac measure supported at $x$.
We note that if $h_{L}(x) = + \infty$, then 
$A^{2}(x,(A + L)_{x},h_{A}\cdot h_{L},\delta_{x}) = 0$.  
For every $x\in X$ we set 
\[
\hat{K}_{m}^{A}(x) := \sup \{ \mid\sigma\mid^{\frac{2}{m}}(x)\mid 
\sigma \in \Gamma (X,{\cal O}_{X}(A + mK_{X})), 
\mid\!\!\int_{X}h_{A}^{\frac{1}{m}}\cdot(\sigma\wedge\bar{\sigma})^{\frac{1}{m}}\!\mid = 1\}.
\]
Here $\mid\sigma\mid^{\frac{2}{m}}$ is not a function on $X$, but the supremum
is takan as a section of the real line bundle $\mid\!A\!\mid^{\frac{2}{m}}\otimes
\mid\!K_{X}\!\mid^{2}$ in the obvious manner\footnote{We have abused the notations $\mid\!\!A\!\!\mid$, $\mid\!\!K_{X}\!\!\mid$ here.  These notations are similar to 
the notations of corresponding linear systems.  But I think there is no fear of confusion. }.  
Then $h_{A}^{\frac{1}{m}}\cdot\hat{K}_{m}^{A}$ is a continuous semipositive $(n,n)$ form on $X$. 
Under the above notations, we have the following theorem.

\begin{theorem}\label{main}
We set 
\[
\hat{K}_{\infty}^{A}:= \limsup_{m\rightarrow\infty}h_{A}^{\frac{1}{m}}\cdot\hat{K}_{m}^{A}
\]
and 
\[
\hat{h}_{can,A} := \mbox{\em the lower envelope of}\,\,\,\hat{K}_{\infty}^{-1}.   
\]
Then $\hat{h}_{can,A}$ is an AZD of $K_{X}$.
And we define 
\[
\hat{h}_{can} : =\mbox{\em the lower envelope of}\,\,\,\inf_{A}\hat{h}_{can,A},
\]
where $\inf$ means the pointwise infimum and $A$ runs all the  
ample line bundles on $X$. 
Then $\hat{h}_{can}$ is  a well defined AZD\footnote{I believe that $\hat{h}_{can,A}$ is already independent of the sufficiently ample line bundle $A$.} depending only on $X$.  
$\square$
\end{theorem}

\begin{definition}{\em {\bf (Supercanonical AZD)}}\label{super} 
We call $\hat{h}_{can}$ in Theorem \ref{main} the supercanonical AZD 
of $K_{X}$.   
And we call the semipositive $(n,n)$ form $\hat{h}_{can}^{-1}$ 
the supercanonical volume form on $X$.  $\square$
\end{definition}
\begin{remark}
Here ``super'' means that corresponding volume form $\hat{h}_{can}^{-1}$
satisfies the inequality :
\[
\hat{h}_{can}^{-1} \geqq h_{can}^{-1}, 
\]
if $X$ has nonnegative Kodaria dimension (cf. Theorem \ref{comparison}). $\square$
\end{remark}
In the statement of Theorem \ref{main}, one may think that $\hat{h}_{can,A}$
may dependent of the choice of the metric $h_{A}$. 
But later we prove that  $\hat{h}_{can, A}$ is independent of the choice of $h_{A}$
(cf. Theorem \ref{uniqueness}).  

\subsection{Variation of the supercanonical AZD $\hat{h}_{can}$}\label{intro}
Let $f : X \longrightarrow S$ be an algebraic fiber space, i.e., 
$X,S$ are smooth projective varieties and  $f$ is a projective morphism with connected fibers. 
Suppose that for a general fiber $X_{s}:= f^{-1}(s)$, $K_{X_{s}}$ is 
pseudoeffective \footnote{This condition is equivalent to the one that 
for some regular fiber $X_{s}$, $K_{X_{s}}$ is pseudoeffective. This is well known. 
For the proof, see  Lemma \ref{ext} below for example.}.   
In this case we may define a singular hermitian metric $\hat{h}_{can}$  
on $K_{X/S}$ similarly as above.   
Then $\hat{h}_{can}$ have a nice properties on $f : X \longrightarrow S$ 
as follows.  
\begin{theorem}\label{family}
Let $f : X \longrightarrow S$ be an algebraic fiber space such that 
for a general fiber $X_{s}$, $K_{X_{s}}$ is pseudoeffective. 
We set $S^{\circ}$ be the maximal nonempty Zariski open subset of $S$ such that $f$ is smooth over $S^{\circ}$ and $X^{\circ} = f^{-1}(S^{\circ})$.
Then there exists a unique singular hermitian metric $\hat{h}_{can}$ on 
$K_{X/S}$ such that 
\begin{enumerate}
\item $\hat{h}_{can}$ has semipositive curvature in the sense of current. 
\item $\hat{h}_{can}\!\mid\!\!X_{s}$ is an AZD of $K_{X_{s}}$ for 
every $s \in S^{\circ}$.
\item There exists the union $F$ of at most countable union of proper subvarieties
of $S$ such that for every $s\in S\,\backslash\, F$, 
\[
\hat{h}_{can}\!\!\mid\!X_{s} \leqq \hat{h}_{can,s}
\]
holds, where $\hat{h}_{can,s}$ denotes the supercanonical AZD of $K_{X_{s}}$. 
\item There exists a subset $G$ of measure $0$
 in $S^{\circ}$, such that for every $s \in S^{\circ}\,\backslash\, G$, 
$\hat{h}_{can}\!\mid\!\!X_{s} = \hat{h}_{can,s}$ holds. 
\end{enumerate} 
 $\square$
\end{theorem}
\begin{remark}
Even for $s\in G$, $\hat{h}_{can}\!\!\mid\!\!X_{s}$ is an AZD of $K_{X_{s}}$
by $2$.  I do not know whether $F$ or $G$ really exists in some cases. $\square$ 
\end{remark}
By Theorem \ref{family} and the $L^{2}$-extension theorem (\cite[p.200, Theorem]{o-t}), we obtain the following corollary immediately. 
\begin{corollary}(\cite{s1,s2,tu3})\label{pg}
Let $f : X \longrightarrow S$ be a smooth projective family over 
a complex manifold $S$. 
Then  plurigenera $P_{m}(X_{s}) := \dim H^{0}(X_{s},{\cal O}_{X_{s}}(mK_{X_{s}}))$ is a locally constant function on $S$ $\square$ 
\end{corollary}
The following corollary is immediate consequence of 
Theorem \ref{family}, since the supercanonical AZD is 
always has minimal singularities (cf. Definition \ref{minAZD} and 
Remark \ref{minimality}). 
\begin{corollary}\label{additive}
Let $f : X \longrightarrow Y$ be an algebraic fiber space. 
Suppose that $K_{X}$ and $K_{Y}$ are pseudoeffective. 
Let $\hat{h}_{can}$ be the canonical singular hermitian metric on $K_{X/Y}$ constructed 
as in Theorem \ref{family}. 
Let $\hat{h}_{can,X},\hat{h}_{can,Y}$ be the supercanonical AZD's of 
$K_{X}$ and $K_{Y}$ respectively. 
Then there exists a positive constant $C$ such that 
\[
\hat{h}_{can,X} \leqq C\cdot \hat{h}_{can}\cdot f^{*}\hat{h}_{can,Y}
\]
holds on $X$. $\square$
\end{corollary}
Cororally \ref{additive} is very close to Iitaka's conjecture which 
asserts  that
\[
\mbox{Kod}(X) \geqq \mbox{Kod}(Y) + \mbox{Kod}(F)
\]
holds for any  algebraic fiber space 
$f : X \longrightarrow Y$, where $F$ is a general fiber of 
$f : X \longrightarrow Y$ and $\mbox{Kod}(M)$ denotes the Kodaira dimension
of a compact complex manifold $M$. \vspace{3mm} \\

In this paper all the varieties are defined over $\mathbb{C}$.  
And we frequently use the classical result  that the supremum of a family of plurisubharmonic functions locally uniformly bounded from above is again plurisubharmonic, if we take 
the uppersemicontinuous envelope of the supremum (\cite[p.26, Theorem 5]{l}). 
For simpliciy, we denote the upper(resp. lower)semicontinuous envelope
simply by the upper(resp. lower) envelope.  We note that this adjustment occurs only on the set of measure $0$. In this paper all the singular 
hermitian metrics are supposed to be lowersemicontinuous. 
	
There are other applications of the supercanonical AZD.  Also it is immediate 
to generalize it to the log category and another generalization involving hermitian
 line bundles with semipositive curvature is also possible.   
These will be discussed in the forthcoming papers.

\section{Proof of Theorem \ref{main}}
In this section we shall prove Theorem \ref{main}. 
We shall use the same notations as in Section \ref{construction}. 
The upper estimate of $\hat{K}_{m}^{A}$ is almost the same as in \cite{tu9}, but the lower estimate of $\hat{K}_{m}^{A}$ requires the $L^{2}$ extension theorem
(\cite{o-t,o}).

\subsection{Upper estimate of $\hat{K}_{m}^{A}$}\label{up}
Let $X$ be as in Theorem \ref{main} and let $n$ denote $\dim X$ 
and let $x\in X$ be an arbitrary point.  
Let $(U,z_{1},\cdots ,z_{n})$ be a coordinate neighbourhood of 
$X$ which is biholomorphic to the unit open polydisk $\Delta^{n}$ 
such that $z_{1}(x)= \cdots = z_{n}(x) = 0$.  

Let $\sigma \in \Gamma (X,{\cal O}_{X}(A + mK_{X}))$. 
Taking $U$ sufficiently small, we may assume that $(z_{1},\cdots ,z_{n})$
is a holomorphic local coodinate on a neighbourhood of the closure of $U$ and 
there exists 
a local holomorphic frame $\mbox{\bf e}_{A}$ of $A$ on a neighbourhood
of the clousure of $U$. 
Then there exists a bounded holomorphic function $f_{U}$ on $U$ such that 
\[
\sigma = f_{U}\cdot \mbox{\bf e}_{A}\cdot (dz_{1}\wedge\cdots \wedge dz_{n})^{\otimes m}\]
holds.  
Suppose that 
\[
\mid\!\!\int_{X}h_{A}^{\frac{1}{m}}\cdot(\sigma\wedge\bar{\sigma})^{\frac{1}{m}}\!\mid = 1
\]
holds.   
Then we see that 
\begin{eqnarray*}
\int_{U}\mid f_{U}(z)\mid^{\frac{2}{m}}d\mu (z) 
& \leqq &
(\inf_{U} h_{A}(\mbox{\bf e}_{A},\mbox{\bf e}_{A}))^{-\frac{1}{m}}\cdot \int_{U}h_{A}(\mbox{\bf e}_{A}\,\mbox{\bf e}_{A})^{\frac{1}{m}}\mid f_{U}\mid^{2}d\mu (z) \\
& \leqq & (\inf_{U} h_{A}(\mbox{\bf e}_{A},\mbox{\bf e}_{A}))^{-\frac{1}{m}}
\end{eqnarray*}
hold, where $d\mu (z)$ denotes the standard Lebesgue measure on the coordinate. Hence by the submeanvalue property of plurisubharmonic functions, 
\[
h_{A}^{\frac{1}{m}}\cdot \mid\sigma \mid^{\frac{2}{m}}(x)
\leqq \{\frac{h_{A}(\mbox{\bf e}_{A},\mbox{\bf e}_{A})(x)}{\inf_{U} h_{A}(\mbox{\bf e}_{A},\mbox{\bf e}_{A})}\}^{\frac{1}{m}}\cdot\pi^{-n}\cdot \mid\!dz_{1}\wedge\cdots\wedge dz_{n}\!\mid^{2}(x)
\]
holds. 
Let us fix a $C^{\infty}$ volume form $dV$ on $X$. 
Since $X$ is compact and every line bundle on a contractible Stein manifold is trivial, we have the following lemma.
\begin{lemma}\label{upper}
 There exists a positive constant $C$ 
independent of the line bundle $A$ and the $C^{\infty}$ metric $h_{A}$ such that 
\[
\limsup_{m\rightarrow\infty}h_{A}^{\frac{1}{m}}\cdot\hat{K}_{m}^{A} \leqq C\cdot dV
\]
holds on $X$.  
$\square$ 
\end{lemma}
\subsection{Lower estimate of $\hat{K}_{m}^{A}$}\label{low}
Let $h_{X}$ be any $C^{\infty}$ hermitian metric on $K_{X}$. 
Let $h_{0}$ be an AZD of $K_{X}$ defined by the lower envelope of : 
\[
 \inf \{ h(x) \mid \mbox{$h$ is a singular hermitian metric on $K_{X}$ with  $\Theta_{h}\geqq 0$,$h \geqq h_{X}$}\}.
\]
Then by the classical theorem of Lelong (\cite[p.26, Theorem 5]{l}) it is easy to verify that $h_{0}$ is an AZD of $K_{X}$ (cf. \cite[Theorem 1.5]{d-p-s}).
$h_{0}$ is of  minimal singularities in the following sense. 

\begin{definition}\label{minAZD}
Let $L$ be a pseudoeffective line bundle on a smooth projective variety $X$.
An AZD $h$ on $L$  is said to be an {\bf AZD of minimal singularities}, if 
for any AZD $h^{\prime}$ on $L$, there exists a positive constant $C$ such that
\[
h \leqq   C \cdot h^{\prime}
\]
holds. $\square$ 
\end{definition}
Let us compare $h_{0}$ and $\hat{h}_{can}$.
  
By the $L^{2}$-extension theorem (\cite{o}), we have the following lemma. 
\begin{lemma}\label{l2ext}
There exists a positive constant $C$ independent of $m$ such that 
\[
K(A + mK_{X}, h_{A}\cdot h_{0}^{m-1}) \geqq C\cdot (h_{A}\cdot h_{0}^{m})^{-1}
\]
holds, where $K(A + mK_{X}, h_{A}\cdot h_{0}^{m-1})$ is the (diagonal part of)
Bergman kernel of $A+ mK_{X}$ with respect to the $L^{2}$-inner product:
\[
(\sigma ,\sigma^{\prime}):= (\sqrt{-1})^{n^{2}}\int_{X}\sigma\wedge\bar{\sigma}^{\prime}
\cdot h_{A}\cdot h_{0}^{m-1}, 
\]
where we have considered $\sigma,\sigma^{\prime}$ as  $A+(m-1)K_{X}$ valued
canonical forms. $\square$
\end{lemma}
{\bf Proof of Lemma \ref{l2ext}}.
By the extremal property of the Bergman kernel (see for example \cite[p.46, Proposition 1.4.16]{kr}) we have that
\begin{equation}\label{extremal}
K(A + mK_{X}, h_{A}\cdot h_{0}^{m-1})(x) = \sup \{\mid\!\sigma (x)\!\mid^{2}
\mid\sigma\in\Gamma (X,{\cal O}_{X}(A + mK_{X})\otimes {\cal I}(h_{0}^{m-1})), 
\parallel\!\sigma\!\parallel = 1\}, 
\end{equation}
holds for every $x\in X$, where $\parallel\!\sigma\!\parallel = (\sigma,\sigma)^{\frac{1}{2}}$.  
Let $x$ be a point such that $h_{0}$ is not $+\infty$ at $x$.
Let $dV$ be an arbitrary $C^{\infty}$ volume form on $X$ as in Section \ref{construction}.  
Then by the $L^{2}$-extension theorem (\cite{o,o-t}) and 
the sufficiently ampleness of $A$ (see Section \ref{construction}),
we  may extend  any  $\tau_{x} \in (A + mK_{X})_{x}$ with $h_{A}\cdot h_{0}^{m-1}\cdot dV^{-1}
(\tau_{x},\tau_{x}) = 1$ to a global section 
$\tau \in  \Gamma (X,{\cal O}_{X}(A + mK_{X})\otimes {\cal I}(h_{0}^{m-1}))$
such that 
\[
\parallel\!\tau\!\parallel \leqq C_{0},
\]
where $C_{0}$ is a positive constant independent of $x$ and $m$.
Let $C_{1}$ be a positive constant such that 
\[
h_{0} \geqq C_{1}\cdot dV^{-1}
\]
holds on $X$.  
By (\ref{extremal}), we obtain the lemma by 
taking $C = C_{0}^{-1}\cdot C_{1}$. $\square$\vspace{3mm} \\

Let $\sigma \in \Gamma (X,{\cal O}_{X}(A + mK_{X})\otimes {\cal I}(h_{0}^{m-1}))$ such that 
\[
(\sqrt{-1})^{n^{2}}\cdot \int_{X}\sigma\wedge\bar{\sigma}\cdot h_{A}\cdot h_{0}^{m-1}
= 1
\]
and 
\[
\mid\sigma\mid^{2}(x) = K(A + mK_{X}, h_{A}\cdot h_{0}^{m-1})(x)
\]
hold, i.e., $\sigma$ is a peak section at $x$. 
Then by the H\"{o}lder inequality we have that  
\begin{eqnarray*}
\mid\!\!\int_{X}h_{A}^{\frac{1}{m}}\cdot (\sigma\wedge\bar{\sigma})^{\frac{1}{m}}\!\mid
& \leqq &(\int_{X}h_{A}\cdot h_{0}^{m}\cdot \mid\sigma\mid^{2}\cdot h_{0}^{-1})^{\frac{1}{m}}
\cdot (\int_{X}h_{0}^{-1})^{\frac{m-1}{m}} \nonumber\\
& \leqq & (\int_{X}h_{0}^{-1})^{\frac{m-1}{m}}
\end{eqnarray*}
hold. 
Hence we have the inequality:  
\begin{equation}\label{b}
\hat{K}_{m}^{A}(x) \geqq  K(A + mK_{X}, h_{A}\cdot h_{0}^{m-1})(x)^{\frac{1}{m}}
\cdot (\int_{X}h_{0}^{-1})^{-\frac{m-1}{m}}
\end{equation}
holds.  Now we shall consider the limit 
\[
\limsup_{m\rightarrow\infty}h_{A}^{\frac{1}{m}}\cdot K(A + mK_{X}, h_{A}\cdot h_{0}^{m-1})^{\frac{1}{m}}.
\]
Let us recall the following result.
\begin{lemma}(\cite[p.376, Proposition 3.1]{dem})\label{dem}
\[
\limsup_{m\rightarrow\infty}h_{A}^{\frac{1}{m}}\cdot K(A + mK_{X}, h_{A}\cdot h_{0}^{m-1})^{\frac{1}{m}}= h_{0}^{-1}
\]
holds. $\square$
\end{lemma}
\begin{remark}
In (\cite[p.376, Proposition 3.1]{dem}, Demailly only considered the local 
version of Lemma \ref{dem}.  But the same proof works in our case 
by the sufficiently ampleness of $A$.  This kind of localization principle 
for Bergman kernels is quite standard. $\square$ 
\end{remark}
In fact the $L^{2}$-extension theorem (\cite{o-t,o}) implies the 
inequality 
\[
\limsup_{m\rightarrow\infty}h_{A}^{\frac{1}{m}}\cdot K(A + mK_{X}, h_{A}\cdot h_{0}^{m-1})^{\frac{1}{m}}
\geqq  h_{0}^{-1}
\]
and the converse inequality is elementary. See \cite{dem} for details 
and applications. 
Hence letting $m$ tend to infinity in (\ref{b}), by Lemma \ref{dem}, we have the following lemma.

\begin{lemma}\label{lower}
\[
\limsup_{m\rightarrow\infty}h_{A}^{\frac{1}{m}}\cdot \hat{K}_{m}^{A} \geqq (\int_{X}h_{0}^{-1})^{-1}\cdot h_{0}^{-1}
\]
holds. $\square$ 
\end{lemma}

By  Lemmas \ref{upper} and \ref{lower}, we see that 
\[
\hat{K}_{\infty}^{A} := \limsup_{m\rightarrow\infty}h_{A}^{\frac{1}{m}}\cdot\hat{K}_{m}^{A}
\]
exists as a bounded semipositive $(n,n)$ form on $X$.
We set  
\[
\hat{h}_{can,A} := \mbox{the lower envelope of}\,\,\,(K_{\infty}^{A})^{-1}. 
\]

\subsection{Independence of $\hat{h}_{can,A}$ from $h_{A}$}

In the above construction, $\hat{h}_{can,A}$ depends on 
the choice of the $C^{\infty}$ hermitian metric $h_{A}$ apriori. 
But actually $\hat{h}_{can,A}$ is independent of the choice of $h_{A}$. 

Let $h_{A}^{\prime}$ be another $C^{\infty}$-hermitian metric on $A$. 
We define 
\[
(\hat{K}_{m}^{A})^{\prime}
:= \sup \{ \mid\sigma\mid^{\frac{2}{m}};\, 
\sigma \in \Gamma (X,{\cal O}_{X}(A + mK_{X})), \mid\!\!\int_{X}(h_{A}^{\prime})^{\frac{1}{m}}\cdot (\sigma\wedge\bar{\sigma})^{\frac{1}{m}}\!\!\mid = 1\}.
\]
We note that the ratio $h_{A}/h^{\prime}_{A}$ is a positive $C^{\infty}$-function 
on $X$ and 
\[
\lim_{m\rightarrow\infty}(\frac{h_{A}}{h^{\prime}_{A}})^{\frac{1}{m}} = 1
\]
uniformly on $X$.   
Since the definitions of $\hat{K}_{m}^{A}$ and $(\hat{K}^{A}_{m})^{\prime}$ 
use the extremal properties, we see easily that for every positive number 
$\varepsilon$, there exists a positive integer $N$ such that 
for every $m \geqq N$
\[
(1-\varepsilon )(\hat{K}_{m}^{A})^{\prime}\leqq \hat{K}_{m}^{A} \leqq (1 + \varepsilon )(\hat{K}_{m}^{A})^{\prime}
\]
holds on $X$. 
Hence we obtain the following uniqueness theorem. 
\begin{theorem}\label{uniqueness}
$\hat{K}_{\infty}^{A}= \limsup_{m\rightarrow\infty}h_{A}^{\frac{1}{m}}\cdot\hat{K}_{m}^{A}$ is independent of the choice of the $C^{\infty}$ hermitian metric
$h_{A}$. 
Hence $h_{can,A}$ is independent of the choice of the $C^{\infty}$ hermitian metric $h_{A}$.  $\square$ 
\end{theorem}

\subsection{Completion of the proof of Theorem \ref{main}}

Let $h_{0}$ be an AZD of $K_{X}$ constructed as in Section \ref{up}.
Then by Lemma \ref{lower} we see that 
\[
\hat{h}_{can,A} \leqq (\int_{X}h_{0}^{-1})\cdot h_{0}
\]
holds. 
Hence we see  
\[
{\cal I}(\hat{h}_{can,A}^{m}) \supseteq {\cal I}(h_{0}^{m})
\]
holds for every $m\geqq 1$. This implies that 
\[
H^{0}(X,{\cal O}_{X}(mK_{X})\otimes{\cal I}(h_{0}^{m}))
\subseteq  H^{0}(X,{\cal O}_{X}(mK_{X})\otimes{\cal I}(\hat{h}_{can,A}^{m}))
\subseteq H^{0}(X,{\cal O}_{X}(mK_{X}))
\]
hold, hence 
\[
H^{0}(X,{\cal O}_{X}(mK_{X})\otimes{\cal I}(\hat{h}_{can,A}^{m}))
\simeq  H^{0}(X,{\cal O}_{X}(mK_{X}))
\]
holds for every $m\geqq 1$. 
And by the construction and the classical theorem of Lelong (\cite[p.26, Theorem 5]{l}) 
stated in Section \ref{intro}, $\hat{h}_{can,A}$ has semipositive curvature in the 
sense of current. 
Hence $\hat{h}_{can,A}$ is an AZD of $K_{X}$ and depends only on $X$ and $A$
by Lemma \ref{uniqueness}.

Let us consider 
\[
\hat{K}_{\infty} : = \sup_{A} \hat{K}_{\infty,A}
\]
where $\sup$ means the pointwise supremum and $A$ runs all the 
sufficiently ample line bundle on $X$. 
Then Lemma \ref{upper}, we see that 
$\hat{K}_{\infty}$ is a well defined semipositive $(n,n)$ form on $X$.
We set 
\[
\hat{h}_{can} := \mbox{the lower envelope of}\,\,\, \hat{K}_{\infty}^{-1}. 
\] 
Then by the construction, $\hat{h}_{can} \leqq \hat{h}_{can,A}$ 
for every ample line bundle $A$. 
Since $\hat{h}_{A}$ is an AZD of $K_{X}$, $\hat{h}_{can}$ 
is also an AZD of $K_{X}$ indeed (again by \cite[p.26, Theorem 5]{l}). 
Since $\hat{h}_{can,A}$ depends only on $X$ and $A$, $\hat{h}_{can}$ is uniquely determined by $X$.  
This completes the proof of Theorem \ref{main}. $\square$ 
\begin{remark}\label{minimality} 
As one see Section \ref{low}, we see that $\hat{h}_{can}$ is 
an AZD of $K_{X}$ of minimal singularities (cf. Definition \ref{minAZD}). $\square$
\end{remark}
\subsection{Comparison of $h_{can}$ and $\hat{h}_{can}$}
Suppose that $X$ has nonnegative Kodaira dimension. 
Then by Theorem \ref{canazd}, we can define the canonical AZD 
$h_{can}$ on $K_{X}$.   We shall compare $h_{can}$ and $\hat{h}_{can}$. 

\begin{theorem}\label{comparison}
\[
\hat{h}_{can,A} \leqq h_{can}
\]
holds on $X$.  In particular
\[
\hat{h}_{can} \leqq h_{can}
\]
holds on $X$  $\square$ 
\end{theorem}
{\bf Proof of Theorem \ref{comparison}}.
If $X$ has  negative Kodaira dimension, then the right hand side 
is infinity.  Hence the ineuqality is trivial. 

Suppose that $X$ has nonnegative Kodaira dimension. 
Let $\sigma\in \Gamma (X,{\cal O}_{X}(mK_{X}))$ be an element such that 
\[
\mid\!\!\int_{X}(\sigma\wedge\bar{\sigma})^{\frac{1}{m}}\!\mid = 1
\]
Let $x\in X$ be an arbitrary point on $X$. Since 
 ${\cal O}_{X}(A)$ is globally generated by the definition of $A$, there exists an element 
 $\tau \in \Gamma (X,{\cal O}_{X}(A))$ such that 
 $\tau (x) \neq 0$ and  $h_{A}(\tau ,\tau )\leqq 1$ on $X$. 
Then we see that 
\[
\int_{X}h_{A}(\tau ,\tau )^{\frac{1}{m}}\cdot (\sigma\wedge\bar{\sigma})^{\frac{1}{m}}\leqq 1
\]   
holds. 
This implies that 
\[
\hat{K}_{m}^{A}(x) \geqq \mid\!\tau (x)\!\mid^{\frac{2}{m}}\cdot K_{m}(x)
\]
holds at $x$.   Noting $\tau (x) \neq 0$,letting $m$ tend to infinity, 
we see that 
\[
\hat{K}^{A}_{\infty}(x) \geqq K_{\infty}(x)
\]
holds.   Since $x$ is arbitrary, this completes the proof of Theorem \ref{comparison}.  $\square$  
\begin{remark}
The equality $h_{can} = \hat{h}_{can}$ implies the abundance of $K_{X}$. 
$\square$
\end{remark}

By the same proof we obtain the following comparison theorem 
(without assuming $X$ has nonnegative Kodaira dimension).
\begin{theorem}\label{monotonicity}
Let $A, B$ a sufficiently ample line bundle on $X$. 
Suppose that $B - A$ is globally generated, then 
\[
\hat{h}_{can,B} \leqq \hat{h}_{can,A}
\]
holds.  $\square$ 
\end{theorem}
\begin{remark}\label{any}
Theorem \ref{monotonicity} implies that 
\[
\hat{h}_{can} = \lim_{\ell\rightarrow\infty}\hat{h}_{can,\ell A}
\]
holds for any  ample line bundle $A$ on $X$. $\square$
\end{remark}

\section{Variation of $\hat{h}_{can}$ under projective deformations} 
In this section we shall prove Theorem \ref{family}.
The main ingredient of the proof is the variation of Hodge structure.

\subsection{Construction of $\hat{h}_{can}$ on a family}\label{confam}

Let $f : X \longrightarrow S$ be an algebraic fiber space 
as in Theorem \ref{family}. 

The construction of $\hat{h}_{can}$ can be performed 
simultaeneously on the family as follows.  The same construction works 
for flat projective family with only canonical singularities. 
But for simplicity we shall work on smooth category.

Let $S^{\circ}$  be the maximal nonempty Zariski open subset of $S$ 
such that $f$ is smooth over $S^{\circ}$ and let us set  
$X^{\circ} := f^{-1}(S^{\circ})$. 

Hereafter we shall assume that $\dim S = 1$.  The general case of Theorem \ref{family} easily follows from just by cutting down $S$ to curves. 
 Let $A$ be a sufficiently ample line bundle on $X$  such that for every pseudoeffective  singular hermitian line bundle $(L,h_{L})$, 
${\cal O}_{X}(A + L)\otimes {\cal I}(h_{L})$ and ${\cal O}_{X}(K_{X} + A+ L)
\otimes {\cal I}(h_{L})$ are  globally generated 
and ${\cal O}_{X_{s}}(A + L\!\mid\!\!X_{s})\otimes {\cal I}(h_{L}\!\mid\!\!X_{s})$ and ${\cal O}_{X_{s}}(K_{X_{s}}+ A + L\!\mid\!\!X_{s})\otimes {\cal I}(h_{L}\!\mid\!\!X_{s})$ 
are globally generated for every $s\in S^{\circ}$ as long as $h_{L}\!\!\mid\!X_{s}$ is well defined. 

Let us assume that there exists a smooth member $D$ of $\mid\!2A\!\mid$ 
such that $D$ does not contain any fiber over $S^{\circ}$.  
Let $\sigma_{D}$ a holomorphic section of $2A$ with divisor $D$. 
We consider the singular hermitian metric
\[
h_{A}:= \frac{1}{\mid\sigma_{D}\mid}
\]
on $A$.     
We set 
\[
E_{m} := f_{*}{\cal O}_{X}(A + mK_{X/S}). 
\]
Since we have assumed that $\dim S = 1$, $E_{m}$ is a vector bundle for every 
$m\geqq 1$.  We denote the fiber of the vector bundle over $s\in S$ by 
$E_{m,s}$. 
Then we shall define the sequence of $\frac{1}{m}A$-valued relative volume forms by  
\[
\hat{K}^{A}_{m,s} := \sup \{ \mid\!\sigma\!\mid^{\frac{2}{m}} ;
\sigma\in E_{m,s}, 
\mid\!\!\int_{X_{s}}h_{A}^{\frac{1}{m}}\cdot(\sigma\wedge\bar{\sigma})^{\frac{1}{m}}\!\!\mid = 1\}
\]
for every $s\in S^{\circ}$.   
This fiberwise construction is different from that in Section \ref{construction}  in the following  two points :
\begin{enumerate}
\item We use the singular metric $h_{A}\!\!\mid\!\!X_{s}$ instead of a $C^{\infty}$ hermitian metric on $A\!\mid\!X_{s}$. 
\item We use $E_{m,s}$ instead of $\Gamma (X_{s},{\cal O}_{X_{s}}(A\!\!\mid\!\!X_{s}
+ mK_{X_{s}}))$. 
\end{enumerate}
We note that the 2nd difference occurs only over at most countable union of 
proper analytic subsets in $S^{\circ}$. 
Since $h_{A}$ is singular, at some point $s\in S^{\circ}$ and for some positive 
integer $m_{0}$, $\hat{K}^{A}_{m_{0},s}$ might be identically $0$ on $X_{s}$.
But for any $s\in S^{\circ}$ we find a positive integer $m_{0}$ such that 
for every $m\geqq m_{0}$, we have 
${\cal I}(h_{A}^{\frac{1}{m}}\!\!\mid\!X_{s}) = {\cal O}_{X_{s}}$ holds for every $m\geqq m_{0}$. 
Hence even in this case we see that $\hat{K}^{A}_{m,s}$ is not identically $0$ for every sufficiently large $m$.

We define the relative $\mid\!A\!\mid^{\frac{2}{m}}$ valued volume form $\hat{K}^{A}_{m}$ by  
\[
\hat{K}_{m}^{A}\!\!\mid\!\!X_{s} := \hat{K}^{A}_{m,s} (s\in S)
\]
and a relative volume form $\hat{K}^{A}_{\infty}$ by 
\[
\hat{K}_{\infty}^{A}\!\!\mid\!\!X_{s} := \limsup_{m\rightarrow\infty}h_{A}^{\frac{1}{m}}\cdot\hat{K}^{A}_{m,s} (s\in S).
\]
Of course the above  construction of $\hat{K}_{m,s}^{A}(s\in S^{\circ})$ works also for $C^{\infty}$ hermitian metric instead of the singular $h_{A}$ as above.
The reason why we use the singular $h_{A}$ is that  we shall use the variation of Hodge structure to prove the plurisubharmonic variation of $\log \hat{K}_{m,s}^{A}$.We may use a $C^{\infty}$ metric with strictly positive curvature on $A$, 
instead of the singular $h_{A}$ as above, if we use the plurisubharmonicity properties
of Bergman kernels (\cite[Theorem 1.2]{be}) instead of Theorem \ref{fujita}.
See Theorem \ref{branch} below.

We define  singular hermitian metrics on 
$A + mK_{X/S}$ by 
\[
\hat{h}_{m,A}:= \mbox{the lower envelope of}\,\,\,(\hat{K}_{m}^{A})^{-1}.
\]

Let us fix a $C^{\infty}$ hermitian metric $h_{A,0}$ on $A$ and 
we set 
\[
\hat{h}_{can,A} := \mbox{the lower envelope of}\,\,\,\liminf_{m\rightarrow\infty}
h_{A,0}^{-\frac{1}{m}}\cdot\hat{h}_{m,A}.
\]
Cleary  $\hat{h}_{can,A}$ does not depend on the choice of $h_{A,0}$ (in this sense, the presence of $h_{A,0}$ is rather auxilary). 
Then we define  
\[
\hat{h}_{can} := \mbox{the lower envelope of}\,\,\, \inf_{A}\hat{h}_{can,A},
\]
where $A$ runs all the ample line bundle on $X$. 
At this moment, $\hat{h}_{can}$ is defined  only on $K_{X/S}\!\!\mid\!\!X^{\circ}$.
The extension of $\hat{h}_{can}$ to the singular hermitian metric on the whole $K_{X/S}$ will be discussed later.
 
\subsection{Semipositivity of the curvature current of $\hat{h}_{m,A}$}\label{semipos}

To prove the semipositivity of the curvature of $\hat{h}_{m,A}$, the following 
theorem is essential.

\begin{theorem}(\cite[p.174,Theorem 1.1]{ka3} see also \cite{f,ka1})\label{fujita}
$\phi : M \longrightarrow C$ be a projective morphism with connected fibers 
from a smooth projective variety $M$ onto a smooth curve $C$. 
Let $K_{M/C}$ be the relative canonical bundle. 
We set $F : = \phi_{*}{\cal O}_{M}(K_{M/C}))$ and let $C^{\circ}$ denote 
the nonempty maximal Zariski open subset of $C$ such that $\phi$ is smooth over 
$C^{\circ}$. 
Let $h_{M/C}$ be the hermitian metric on $F\mid C^{\circ}$ by  
\[
h_{M/C}(\sigma ,\sigma^{\prime} ) := (\sqrt{-1})^{n^{2}}\int_{M_{t}}\sigma\wedge\overline{\sigma^{\prime}},
\]
where $n = \dim M - 1$. 
Let $\pi : \mathbb{P}(F^{*})\longrightarrow C$ be the projective bundle
associated with $F^{*}$ and 
Let $L \longrightarrow \mathbb{P}(F^{*})$ be the tautological line bundle.  
Let $h_{L}$ denote the hermitian metric on $L\mid\pi^{-1}(S^{\circ})$ induced by $h_{M/C}$. 

Then $h_{L}$ has semipositive curvature on $\pi^{-1}(S^{\circ})$ and $h_{L}$ extends to the singular hermitian metric on $L$ with semipositive
curvature current. $\square$ 
\end{theorem}

We define the pseudonorm $\parallel\!\sigma\!\parallel_{\frac{1}{m}}$ 
of $\sigma \in E_{m,s}$ by 
\[
\parallel\!\sigma\!\parallel_{\frac{1}{m}}:= \mid\!\!\int_{X_{s}}h_{A}^{\frac{1}{m}}\cdot (\sigma\wedge\bar{\sigma})^{\frac{1}{m}}\!\mid^{\frac{m}{2}}.
\]

We set $E_{m} = f_{*}{\cal O}_{X}(A + mK_{X/S})$ 
and let $L_{m}$ be the tautological line bundle 
on $\mathbb{P}(E_{m}^{*})$, where $E_{m}^{*}$ denotes the dual of 
$E_{m}$. 
By Theorem \ref{fujita} and the branched covering trick, we obtain 
the following essential lemma.   
\begin{lemma}(\cite[p.63, Lemma 7 and p.64, Lemma 8]{ka1})\label{subh}
Let \\ $\sigma \in \Gamma (X,{\cal O}_{X}(A + mK_{X/S}))$.
Then  $\parallel\,\,\,\,\,\,\parallel$ 
defines a singular hermitian metric with semipositive curvature 
on $L_{m}$. 
$\square$
\end{lemma}  
{\bf Proof of Lemma \ref{subh}.}
If  there were no $A$,  the lemma is  completely the same as 
\cite[p.63, Lemma 7 and p.64, Lemma 8]{ka1}.  
In our case, we use the Kawamata's trick to reduce the logarithmic case 
to the non logarithmic case.   Since this trick has been used repeatedly 
by Kawamata himself (see \cite{ka2,ka3} for example), the following argument has no originality. 
We consider the multivalued relative log canonical form 
\[
(\frac{\sigma}{\sqrt{\sigma_{D}}})^{\frac{1}{m}}.  
\]
Then there exists a finite cyclic covering 
\[
\mu : Y \longrightarrow X
\]
such that $\mu^{*}(\frac{\sigma}{\sqrt{\sigma_{D}}})^{\frac{1}{m}}$
is a (single valued) relative canonical form on $Y$ \footnote{If we use 
a $C^{\infty}$ hermitian metric instead of the above $h_{A}$, we also 
construct a cyclic covering $\mu : Y \longrightarrow X$ such that $\frac{1}{m}\mu^{*}L$ is a genuine line bundle on $Y$ and 
$\mu^{*}\sigma^{\frac{1}{m}}$ is a $\frac{1}{m}\mu^{*}L$-valued 
canonical form on $Y$.}.
Here the branch locus of $\mu$ may be  much larger than the union of 
$D \cup (\sigma )$.   But it does not matter. 
The branch covering is used  to reduce the log canonical case 
to the canonical case. 
Let $\pi :  \tilde{Y}\longrightarrow Y$ be an equivariant resolution of singularities and 
let 
\[
\tilde{f} : \tilde{Y} \longrightarrow S
\]
be the resulting family. 
We shall denote the composition $\mu\circ\pi : \tilde{Y}\longrightarrow X$
by $\tilde{\mu}$. 
Let $U$ be a Zariski open subset of $S_{\sigma}$ such that 
$\tilde{f}$ is smooth.
 We note 
that the Galois group action is isometric on $\tilde{f}_{*}{\cal O}_{\tilde{Y}}(K_{\tilde{Y}/S})$
with respect to the natural $L^{2}$-inner product on 
$\tilde{f}_{*}{\cal O}_{\tilde{Y}}(K_{\tilde{Y}/S})$.  
Therefore by Theorem \ref{fujita}, we see that $\parallel\,\,\,\,\,\,\parallel$  defines a singular hermitian metric on $L_{m}$ with semipositive curvature on a nonempty Zariski open  
of $\mathbb{P}(E_{m}^{*})$. 
Again by  Theorem \ref{fujita} the singular hermitian metric extends to 
the whole $\mathbb{P}(E^{*}_{m})$ preserving semipositive curvature property. \vspace{3mm}  $\square$  \\
We also present an alternative proof indicated by Bo Berndtsson
at the workshop at MSRI in April, 2007. \vspace{3mm}\\ 
{\bf Alternative proof of Lemma \ref{subh}}(cf. \cite[Section 6]{b-p}).\\ 
We use the eqality 
\[
\mid\sigma\mid^{\frac{2}{m}} = \frac{\mid\sigma\mid^{2}}{\mid\sigma\mid^{2\frac{m-1}{m}}}
\]
and view 
\[
h_{A}^{\frac{1}{m}}\cdot\frac{1}{\mid\sigma\mid^{2\frac{m-1}{m}}}
\]
as a singular hermitian metric on $(m-1)K_{X/S}+A$. 
Then by \cite[Therem 5.4]{tu8} or \cite{b-p}, we see that
\[
\mid\int_{X_{s}}h_{A}^{\frac{1}{m}}\cdot(\sigma\wedge\bar{\sigma})^{\frac{1}{m}}\mid^{\frac{2}{m}}
\]
defines a singular hermitian metric with semipositive curvature current on $L_{m}$. 
The rest of the proof is identical as the previous one.  $\square$
\begin{remark}\label{rmsubh}
The metric $h_{A}$ can be replaced by a $C^{\infty}$-hermitian metric 
with semipositive  curvature in the second proof.  
$\square$ 
\end{remark}
\begin{corollary}(see also \cite[Section 6]{b-p}) \label{poscurv}
The curvature $\Theta_{\hat{h}_{m,A}} = \sqrt{-1}\partial\bar{\partial}\log \hat{K}_{m,A}$ is semipositive everywhere on $X^{\circ}$.
 $\square$ 
\end{corollary}
{\bf Proof}. Let $x \in X_{s}(s\in S^{\circ})$ and let $\Omega$ be a holomorphic local generatorof $K_{X/S}$ and let $\mbox{\bf e}_{A}$ be a holomorphic local generator of $A$ 
on a neighbourhood $U$ of $x$ in $X^{\circ}$.  
Viewing $\xi(y):= (\mbox{\bf e}_{A}^{-1}\cdot \Omega^{-m})(y)$ as an element of the dual 
of $E_{m,f(y)}$ by $\sigma \in E_{m,f(y)} \mapsto 
\sigma (y)\cdot (\mbox{\bf e}_{A}^{-1}\cdot \Omega^{-m})(y) (y\in U)$,  
\[
\log (\hat{K}_{m,A}(y)\cdot \mid \mbox{\bf e}_{A}\mid^{-\frac{2}{m}}\cdot \mid\Omega\mid^{-2}(y))
\,\,\, (y\in U) 
\]
is plurisubharmonic function on $U$, since
\[
\mid\xi(y)\mid^{\frac{2}{m}}\cdot \hat{K}_{m,A}(y) = 
\sup \{\frac{\mid\xi(y)\cdot \sigma(y)\mid^{\frac{2}{m}}}{\parallel[\sigma]_{[\xi(y)]}\parallel^{\frac{2}{m}}} ; \,\,
\sigma\in E_{m,f(y)}, [\sigma]_{[\xi(y)]}\neq 0 \} 
\]
holds, where $[\sigma]_{[\xi(y)]}$ denotes the class of $\sigma\in E_{m,f(y)}$
in the fiber $L_{m,[\xi(y)]}$ 
at $[\xi(y)]\in \mathbb{P}(E_{m}^{*})$. 
$\square$ \vspace{5mm}\\ 
Now let us consider the behavior of $\hat{h}_{m,A}$ along $X \backslash X^{\circ}$.
Since the problem is local, we may and do assume $S$ is a unit open
disk $\Delta$ in $\mathbb{C}$ for the time being. 
For every local holomorphic section $\sigma$ of $E_{m}$  the function 
\[
\mid\!\!\int_{X_{s}}h_{A}^{\frac{1}{m}}\cdot (\sigma\wedge\bar{\sigma})^{\frac{1}{m}}\!\!\mid
\] 
is of algebraic growth along $S\,\backslash\, S^{\circ}$.    
More precisely for $s_{0} \in S\,\backslash\, S^{\circ}$ as in \cite[p.59 and p. 66]{ka1} there exist positive numbers $C,\alpha ,\beta$ such that 
\begin{equation}\label{eq}
\mid\!\!\int_{X_{s}}h_{A}^{\frac{1}{m}}\cdot (\sigma\wedge\bar{\sigma})^{\frac{1}{m}}\!\mid 
\leqq C\cdot \mid\!\!s - s_{0}\!\!\mid^{-\alpha}\cdot\mid\!\log (\!s-s_{0}\!)\!\mid^{\beta}
\end{equation}
holds. 
Moreover as  \cite[p.66]{ka1} for a nonvanishing 
holomorphic section $\sigma$ of $E_{m}$ around $p\in S\,\backslash\, S^{\circ}$,  the pseudonorm 
\[
\parallel\!\sigma\!\parallel_{\frac{1}{m}}= \mid\!\!\int_{X_{s}}h_{A}^{\frac{1}{m}}(\sigma\wedge\bar{\sigma})^{\frac{1}{m}}\!\!\mid^{\frac{m}{2}}
\]
has a positive lower bound around every $p\in S$. 
This implies that $\hat{h}_{m,A}$ is bounded from below by a smooth metric 
along the boundary $X\,\backslash\, X^{\circ}$. 
By the above estimate, $\hat{h}_{m,A}$ is  of algebraic growth  along the fiber on $X\,\backslash\, X^{\circ}$ by its definition     
and $\hat{h}_{m,A}$ extends to a singular hermitian metric on $\frac{1}{m}A+K_{X/S}$ with semipositive curvature on the whole $X$.   

Now we set 
\[
\hat{h}_{can,A} := \mbox{the lower envelope of}\,\,\,\liminf_{m\rightarrow\infty}h_{A,0}^{-\frac{1}{m}}\cdot\hat{h}_{m,A}, 
\]
where $h_{A,0}$ be a $C^{\infty}$ metric on $A$ (with strictly positive curvature) as in the last subsetion \footnote{One may use $h_{A}$ instead of $h_{A,0}$ here.  But the corresponding limits may be different along $D$, although 
the difference is negligible by taking the lower envelopes. }.

To extend $\hat{h}_{can,A}$ across $S\,\backslash\, S^{\circ}$, we use the 
following useful lemma.

\begin{lemma}(\cite[Corollary 7.3]{b-t})\label{b-t}
Let $\{ u_{j}\}$ be a sequence of plurisubharmonic functions locally bounded above 
on the bounded open set $\Omega$ in $\mathbb{C}^{m}$. Suppose further 
\[
\limsup_{j\rightarrow\infty}u_{j}
\]
is not identically $-\infty$ on any component of $\Omega$. 
Then there exists a plurisubharmonic function $u$ on $\Omega$
such that the set of points
\[
\{x \in \Omega \mid u(x) \neq (\limsup_{j\rightarrow\infty}u_{j})(x)\}
\]
is pluripolar.
$\square$
\end{lemma}

Since $\hat{h}_{m,A}$ extends to a singular hermitian metric on $\frac{1}{m}A + K_{X/S}$
with semipositive curvature current  on the whole $X$ and 
\[
\hat{h}_{can,A}:= \mbox{the lower envelope of}\,\,\,\liminf_{m\rightarrow\infty}h_{A,0}^{-\frac{1}{m}}\cdot\hat{h}_{m,A}
\]
exists as a singular hermitian metric on $K_{X/S}$ on $X^{\circ} = f^{-1}(S^{\circ})$, we see that $\hat{h}_{can,A}$ extends as a singular hermitian metric 
with semipositive curvature current on the whole $X$ by Lemma \ref{b-t}.

Repeating the same argument we see that $\hat{h}_{can}$ is a well defined 
singular hermitian metric with semipositive curvature current on 
$K_{X/S}\!\mid\!X^{\circ}$ and it extends to a singular hermitian metric 
on $K_{X/S}$ with semipositive curvature current on the whole $X$.

\subsection{Uniqueness of $\hat{h}_{can,A}$ for singular $h_{A}$'s}
In the above construction, we use a singular hermitian metric $h_{A}$ on 
$A$ instead of a $C^{\infty}$ hermitian metric. 
We note that $h_{A}$ is singular along the divisor $D$. 
Hence the resulting metric may be a little bit different from the original construction apriori. 
But actually Theorem \ref{uniqueness} still holds. 
Our metric $h_{A}$ is defined as 
\[
h_{A} = \frac{1}{\mid\!\sigma_{D}\!\!\mid}
\]
as above.  Let $h_{A}^{\prime}$ be a $C^{\infty}$ hermitian metric on 
$A$.  
Let us fix an arbitrary point  $s\in S^{\circ}$. 
Let us fix a K\"{a}hler metric on $X$ and let $U_{\varepsilon}$ be 
the $\varepsilon$ neighbourhood of $D$ with respect to the metric. 
 By the upper estimate Lemma \ref{upper}, we see that 
although $h_{A}$ is singular along $D$, there exists a positive integer $m_{0}$ and a positive constnat $C$ 
depending only on $s$ such that for every $m\geqq m_{0}$ and  any $\sigma \in E_{m,s}$ with 
\[
\parallel\sigma\parallel_{\frac{1}{m}}
= \mid\!\!\int_{X_{s}}h_{A}^{\frac{1}{m}}\cdot (\sigma\wedge\bar{\sigma})^{\frac{1}{m}}\!\!\mid^{\frac{m}{2}} = 1,
\]
\[
\mid\!\!\int_{U_{\varepsilon}\cap X_{s}}h_{A}^{\frac{1}{m}}\cdot (\sigma\wedge\bar{\sigma})^{\frac{1}{m}}\!\!\mid \leqq C\cdot \varepsilon
\]
holds. 
This means that there is no mass concentration around the neighbourhood
of $D\cap X_{s}$.  We note that on $X_{s}\,\backslash\, U_{\varepsilon}$
the ratio $(h_{A}/h_{A}^{\prime})^{\frac{1}{m}}$ converges uniformly 
to $1$ as $m$ tends to infinity. 
Hence by the definitions of $\hat{K}_{m,s}^{A}$ and 
$(\hat{K}^{A}_{m,s})^{\prime}$, we see that for every $s\in S^{\circ}$ and 
$\delta > 0$,  there exists a positive integer 
$m_{1}$ such that for every $m \geqq m_{1}$
\[
(1-\delta )(\hat{K}_{m,s}^{A})^{\prime} \leqq \hat{K}_{m,s}^{A} \leqq (1 + \delta )(\hat{K}_{m,s}^{A})^{\prime}
\]
holds on $X_{s}$. 
Hence  we have the following lemma. 
\begin{lemma}\label{uniqueness2}
$\hat{K}_{\infty,s}^{A}$ is same as the one 
defined by a $C^{\infty}$ hermitian metric on $A$ for every $s\in S^{\circ}$. $\square$
\end{lemma} 

\subsection{Case $\dim S > 1$}
In Sections \ref{confam},\ref{semipos}, 
we have assumed that $\dim S = 1$.  
In the case of $\dim S > 1$ the same proof works similarly.  
But there are several minor differences. 

First there may not exist $D\in \mid 2A\mid$ 
which does not contain any fibers, hence the restriction of $h_{A}$ 
may not be well defined on some fibers in this case. 
But this can be taken care by Lemma \ref{uniqueness2}. 
Namely $\hat{h}_{can}$ is independent of the choice of $D$.  
Hence replacing $h_{A}$ by a $C^{\infty}$ hermitian metric,  
we see that $\hat{K}_{\infty}^{A}$ is defined on all fibers over $S^{\circ}$. 

Second  in this case $E_{m} = f_{*}{\cal O}_{X}(A + mK_{X/S})$ may not be 
locally free on $S^{\circ}$.   If $E_{m.s}$ is not locally free at $s_{0}\in S^{\circ}$, then $\hat{K}_{\infty}^{A}$ may be discontinuous at $s_{0}$.
But 
\[
J:= \{s \in S^{\circ}\mid E_{m} \,\,\,\mbox{is not locally free at $s$ for some $m\geqq 1$}\}
\]
is at most a countable union of proper subvarieties  
of $S^{\circ}$ and 
\[
\hat{h}_{can,A}:= \mbox{the lower envelope of}\,\,\frac{1}{\hat{K}_{\infty}^{A}}
\]
is a well defined singular hermitian metric with semipositive curvature current
on $X^{\circ}$, i.e., the construction is indifferent to the thin set $J$.    
Hence we may construct $\hat{h}_{can}$ on $X^{\circ}$ in this case. 
The extension of $\hat{h}_{can}$ as a singular hermitian metric on $K_{X/S}$
with semipositive curvature current can be accomplished just by slicing $S$ 
by curves.    Hence we complete the proof of the assertion $1$ in Theorem \ref{family}.

\subsection{Completion of the proof of Theorem \ref{family}}

To complete the proof of Theorem \ref{family}, we need to show 
that $\hat{h}_{can}$ defines an AZD for $K_{X_{s}}$ for every $s\in S$. 
To show this fact, we modify the construction of $\hat{K}_{m}^{A}$. 
Here we do not assume  $\dim S = 1$. 

Let us fix $s \in S^{\circ}$ and let $h_{0,s}$ be an AZD constructed as 
in Section \ref{low}.   
Let $U$ be a neighbourhood of $s\in S^{\circ}$ in $S^{\circ}$ which is biholomorphic to an open ball in $\mathbb{C}^{k} (k:= \dim S)$.
By the $L^{2}$-extension theorem (\cite{o-t,o}), we have the following lemma. 
\begin{lemma}\label{ext}
Every element of 
$\Gamma (X_{s},{\cal O}_{X_{s}}(A\mid X_{s} + mK_{X_{s}})
\otimes {\cal I}(h_{0,s}^{m-1}))$ 
extends to an element of  
$\Gamma (f^{-1}(U),{\cal O}_{X}(A + mK_{X}))$
for every positive integer $m$. $\square$ 
\end{lemma}

\noindent{\bf Proof of Lemma \ref{ext}}. 
We prove the lemma by induction on $m$.
If $m=1$, then the $L^{2}$-extension theorem (\cite{o-t,o}) 
implies that every element of 
$\Gamma (X_{s},{\cal O}_{X_{s}}(A + K_{X_{s}}))$ 
extends to an element of $\Gamma (f^{-1}(U),{\cal O}_{X}(A + K_{X}))$.
 Let $\{\sigma_{1,s}^{(m-1)},\cdots ,\sigma_{N(m-1)}^{(m-1)}\}$
be a basis of  $\Gamma (X_{s},{\cal O}_{X_{s}}(A\mid X_{s} + (m-1)K_{X_{s}})
\otimes {\cal I}(\tilde{h}_{0,s}^{m-2}))$ for some $m \geqq 2$. 
Suppose that we have already constructed holomorphic extensions
\[
\{\tilde{\sigma}_{1,s}^{(m-1)},\cdots ,\tilde{\sigma}_{N(m-1),s}^{(m-1)}\}
\subset \Gamma (f^{-1}(U),{\cal O}_{X}(A + (m-1)K_{X}))
\]
of   $\{\sigma_{1,s}^{(m-1)},\cdots ,\sigma_{N(m-1),s}^{(m-1)}\}$
to $f^{-1}(U)$. 
We define the singular hermitian metric $H_{m-1}$ on $(A + (m-1)K_{X})\mid f^{-1}(U)$ by 
\[
H_{m-1} :=
 \frac{1}{\sum_{j=1}^{N(m-1)}\mid\!\tilde{\sigma}^{(m-1)}_{j,s}\!\!\mid^{2}}.
\]
We note that by the choice of $A$, 
${\cal O}_{X_{s}}(A\!\mid\!\!X_{s} + mK_{X_{s}})\otimes {\cal I}(h_{0,s}^{m-1})$ 
is globally generated. 
Hence we see that 
\[
{\cal I}(h_{0,s}^{m})\subseteq {\cal I}(h_{0,s}^{m-1})\subseteq {\cal I}(H_{m-1}\!\!\mid\!X_{s}) 
\]
hold on $X_{s}$.  Apparently $H_{m-1}$ has a semipositive curvature current. 
Hence by the $L^{2}$-extension theorem (\cite[p.200, Theorem]{o-t}), we  may extend every element of
\[ 
\Gamma (X_{s},{\cal O}_{X_{s}}(A + mK_{X_{s}})\otimes {\cal I}(h_{0,s}^{m-1}))
\]
extends to an element of 
\[
\Gamma (f^{-1}(U),{\cal O}_{X}(A + mK_{X})\otimes {\cal I}(H_{m-1})).
\]
This completes the proof of Lemma \ref{ext} by induction. $\square$  \vspace{3mm} \\

Let $h_{A,0}$ be a $C^{\infty}$ hermitian metric on $A$ 
with strictly positive curvature as in the end of the last subsection. 
We define the sequence of $\{ \tilde{K}^{A}_{m,s}\}$  by  
\[
\tilde{K}_{m,s}^{A} := \sup \{ \mid\sigma\mid^{\frac{2}{m}} ;\,
\sigma\in \Gamma (X_{s},{\cal O}_{X_{s}}(A\mid X_{s} + mK_{X_{s}})
\otimes {\cal I}(h_{0,s}^{m-1})),
\mid\!\!\int_{X_{s}}h_{A,0}^{\frac{1}{m}}\cdot(\sigma\wedge\bar{\sigma})^{\frac{1}{m}}\!\!\mid = 1\}.
\]
\noindent By Lemma \ref{ext}, we obtain the following lemma 
immediately. 
\begin{lemma}\label{limit}
\[
\limsup_{m\rightarrow\infty}h_{A,0}^{\frac{1}{m}}\cdot\tilde{K}^{A}_{m,s} \leqq \hat{K}^{A}_{\infty,s}
\]
holds.  $\square$
\end{lemma}
{\bf Proof}. 
We set 
\[
\hat{K}^{A,0}_{m,s} = \sup \{ \mid\sigma\mid^{\frac{2}{m}}
; \,\sigma\in E_{m,s}, \mid\!\!\int_{X_{s}}h_{A,0}^{\frac{1}{m}}\cdot 
(\sigma\wedge \bar{\sigma})^{\frac{1}{m}}\!\!\mid = 1\}.
\]
Then by the definition of $\tilde{K}^{A}_{m,s}$ and Lemma \ref{ext}
we have that 
\begin{equation}\label{3}
\tilde{K}^{A}_{m,s} \leqq \hat{K}^{A,0}_{m,s} 
\end{equation}
holds on $X_{s}$. 
On the other hand by  Lemma \ref{uniqueness2}, we see that  
\begin{equation}\label{4}
\limsup_{m\rightarrow\infty} h_{A,0}^{\frac{1}{m}}\cdot\hat{K}^{A,0}_{m,s}
= \limsup_{m\rightarrow\infty}h_{A,0}^{\frac{1}{m}}\cdot\hat{K}^{A}_{m,s} = \hat{K}_{\infty,s}
\end{equation}
hold. 
Hence combining (\ref{3}) and (\ref{4}), we complete the proof of Lemma \ref{limit}. $\square$ \vspace{3mm}\\

\noindent We set 
\[
\tilde{h}_{m,A,s} := (\tilde{K}^{A}_{m,s})^{-1}. 
\]
We have the following lemma. 
\begin{lemma}
If we define 
\[
\tilde{K}^{A}_{\infty ,s}:= \limsup_{m\rightarrow\infty}h_{A,0}^{\frac{1}{m}}\cdot\tilde{K}^{A}_{m,s}
\]
and 
\[
\tilde{h}_{\infty,A,s} := \mbox{{\em the lower envelope of}}\,\,\tilde{K}_{\infty .A,s}^{-1},
\]
$\tilde{h}_{\infty,A,s}$ is an AZD of $K_{X_{s}}$. $\square$ 
\end{lemma}
{\bf Proof}. 
Let $h_{0,s}$ be an AZD of $K_{X_{s}}$ as above. 
We note that  \\
${\cal O}_{X_{s}}(A\!\!\mid\!\!X_{s} + mK_{X_{s}})\otimes {\cal I}(h_{0,s}^{m-1})$
is globally generated by the definition of $A$. 
Then by the definition of $\tilde{K}_{m,s}^{A}$,  
\[
{\cal I}(h_{0,s}^{m}) \subseteq {\cal I}(\tilde{h}_{m,A,s}^{m})
\]
holds for every $m\geqq 1$.  Hence by repeating the arugument in Section \ref{low},
similar to  Lemma \ref{lower}, we have that 
\[
\tilde{h}_{\infty,A,s} \leqq (\int_{X_{s}}h_{0,s}^{-1})\cdot h_{0,s}
\]
holds.    Hence $\tilde{h}_{\infty,A,s}$ is an AZD of $K_{X_{s}}$. 
$\square$  \vspace{7mm} \\
Since by the construction  and Lemma \ref{uniqueness2}
\[
\hat{h}_{can,s} \leqq \tilde{h}_{\infty,A,s}
\]
holds on $s$, 
we see that $\hat{h}_{can}\!\mid\!X_{s}$ is an AZD of $K_{X_{s}}$.
Since $s\in S^{\circ}$ is arbitrary, we see that 
$\hat{h}_{can}\!\!\mid\!\!X_{s}$ is an AZD of $K_{X_{s}}$ for every $s\in S^{\circ}$.  This completes the proof of the assertion $2$ in Theorem \ref{family}. 
We have already seen that the singular hermitian metric $\hat{h}_{can}$ has semipositive curvature in the sense of current (cf. Section \ref{semipos} expecially Corollary \ref{poscurv}). 
We note that there exists the union $F$ of at most countable union of 
proper subvarieties of $S$ such that for every $s \in S^{\circ}\,\backslash\, F$
\[
E^{(\ell)}_{m,s} = \Gamma (X_{s},{\cal O}_{X}(\ell A + mK_{X_{s}}))
\]
holds for every $\ell , m\geqq 1$. 
Then by the construction and Theorem \ref{monotonicity}(see Remark \ref{any})\footnote{\noindent Theorem \ref{monotonicity} is used because some ample line bundle on the fiber may not extends to an ample line bundle on $X$ in general.} for every $s\in S^{\circ}\,\backslash\,F$, 
\[
\hat{h}_{can}\!\!\mid\!X_{s} \leqq  \hat{h}_{can,s}
\]
holds, where $\hat{h}_{can,s}$ is the supercanonical AZD of $K_{X_{s}}$. 
This completes the proof of the assertion 3 in Theorem \ref{family}. 

We shall define the singular hermitian metric $\hat{H}_{can}$ on 
$K_{X/S}\!\!\mid\!X^{\circ}$ by 
\[
\hat{H}_{can}\!\!\mid\!X_{s}:= \hat{h}_{can,s}\hspace{5mm} (s\in S^{\circ}).
\]
Then by the construction of $\hat{h}_{can}$ there exists a subset $Z$ of measure $0$ in $X^{\circ}$, 
such that 
\[
\hat{H}_{can}\!\!\mid\!X^{\circ}\,\backslash\,Z = \hat{h}_{can}\!\!\mid\!X^{\circ}\,\backslash\,Z
\]
holds.  Let us set
\[
G := \{ s\in S^{\circ}\mid \mbox{$X_{s}\cap Z$ is not of measure $0$ 
in $X_{s}$}\}.
\] 
Then since $Z$ is of measure $0$, $G$ is of measure $0$ in $S^{\circ}$.   
For $s \in S\,\backslash\, G$, by the definition of the supercanonical 
AZD $\hat{h}_{can,s}$ of $K_{X_{s}}$, we see that 
\[
\hat{h}_{can}\!\!\mid\!X_{s} = \hat{h}_{can,s}
\]
holds.
This completes the proof of Theorem \ref{family}. $\square$ 

\begin{remark} As above we have used the singular hermitian metric $h_{A}$ to 
prove Theorem \ref{family} and then go back to the case of a $C^{\infty}$ metric by the uniqueness result (Lemma \ref{uniqueness2}). 
This kind of interaction between singular and smooth metrics have been 
seen in the convergence of the currents associated with random sections of a positive line bundle to the 1-st Chern 
form of the positive line bundle (see \cite{s-z}). 
My first plan of the proof of   Theorem \ref{family} was to use the random sections 
to go to the smooth case from the singular case.  Although I cannot justify it, it seems to be 
interesting to pursue this direction. $\square$ 
\end{remark} 

\section{Appendix}
The following theorem is a generalization of Theorem \ref{fujita}. 
\begin{theorem}\label{branch}
$\phi : M \longrightarrow C$ be a projective morphism with connected fibers 
from a smooth projective variety $M$ onto a smooth curve $C$. 
Let $K_{M/C}$ be the relative canonical bundle.
Let $(L,h_{L})$ be a pseudoeffective singular hermitian line bundle on $M$ Let $m$ be a positive integer. 
We set $F : = \phi_{*}{\cal O}_{M}(mK_{M/C}+L)$ and let $C^{\circ}$ denote 
the nonempty maximal Zariski open subset of $C$ such that $\phi$ is smooth over  $C^{\circ}$. 
Let $\pi : \mathbb{P}(F^{*})\longrightarrow C$ be the projective bundle
associated with $F^{*}$ and 
Let $H \longrightarrow \mathbb{P}(F^{*})$ be the tautological line bundle.  
Let $h_{H}$ denote the singular hermitian metric on $H\mid\pi^{-1}(S^{\circ})$ 
defined by  
\[
h_{H}(\sigma ,\sigma) := \{(\sqrt{-1})^{n^{2}}\int_{M_{t}}
h_{L}^{\frac{1}{m}}\cdot (\sigma\wedge\overline\sigma)^{\frac{1}{m}}\}^{\frac{m}{2}},
\]
where $n = \dim M - 1$. 
Then $h_{H}$ has semipositive curvature on $\pi^{-1}(S^{\circ})$ and $h_{H}$ extends to the singular hermitian metric on $H$ with semipositive
curvature current. $\square$ 
\end{theorem}
{\bf Proof}. The proof is a minor modification  of the proof of 
Lemma \ref{subh}. 
Let $\sigma$ be a local holomorphic section of $H$ on $\pi^{-1}(S^{\circ})$. 
We consider the multivalued $\frac{1}{m}L$-valued canonical form
$\sqrt[m]{\sigma}$ and uniformize it by taking a suitable cyclic Galois covering 
\[
\mu : Y \longrightarrow X
\]
as in Lemma \ref{subh}. 
Then applying  \cite[Theorem 1.2]{be} or \cite[Theorem 5.4]{tu8} (see also \cite{b-p}) on $Y$, as in Lemma \ref{subh}, we see that $h_{H}$ defines a singular hermitian metric 
on the tautological line bundle on $\mathbb{P}(E_{m}^{*})$.  Hence we see that $h_{H}$ has semipositive curvature on  $\pi^{-1}(S^{\circ})$. The extension of $h_{H}$ to the whole $H$ is also follows from \cite[Theorem 5.4]{tu8}. This completes the proof of Theorem \ref{branch}. $\square$  

\small{

}

\noindent Author's address\\
Hajime Tsuji\\
Department of Mathematics\\
Sophia University\\
7-1 Kioicho, Chiyoda-ku 102-8554\\
Japan
\end{document}